\documentclass[a4wide,10pt]{article}%
\usepackage{geometry}                
\usepackage{tikz} 
\usetikzlibrary{positioning}
\usepackage{graphicx}
\usepackage{authblk}
\usepackage{amssymb}
\usepackage{epstopdf}
\usepackage{subfigure}  
\usepackage[sans]{dsfont}
\usepackage[applemac]{inputenc}
\usepackage[english]{babel}
\usepackage{latexsym}
\usepackage{mathrsfs}
\usepackage{amscd}
\usepackage{graphicx}
\usepackage{color}
\usepackage{float}
\usepackage{amsmath}
\usepackage{amsfonts}
\numberwithin{equation}{section}
\usepackage{enumerate}
\usepackage{amsthm}
\usepackage{algorithm}
\usepackage{algorithmic}
\usepackage[numbers,sort]{natbib}
\usepackage[bookmarks=true,colorlinks=true,linkcolor={blue},urlcolor={blue}, citecolor={blue},pdfstartview={XYZ null null 1.22}]{hyperref}%

\theoremstyle{definition}
\newtheorem{alg}{Algorithm}[section]
\newtheorem{remark}{Remark}
\newcommand{\xx}{{\bf x}}

\def\RR{\mathbb R}

\newcommand{\T}{\mathcal T}

\def\be{\begin{equation}}
	\def\ee{\end{equation}}
\def\bea{\begin{eqnarray}}
	\def\eea{\end{eqnarray}}

\title{Localized KBO with genetic dynamics for multi-modal optimization}

\author{ Federica Ferrarese\footnote{		Department of Mathematics and Computer Science, University of Ferrara, e-mail: federica.ferrarese@unife.it}, Claudia Totzeck\footnote{		School of Mathematics and Natural Sciences, University of Wuppertal, e-mail: 	totzeck@uni-wuppertal.de}}

\begin{document}
	\date{}
	\maketitle

	\begin{abstract}
		In this paper, we introduce a novel approach to multi-modal optimization by enhancing the recently developed kinetic-based optimization (KBO) method with genetic dynamics (GKBO). The proposed method targets objective functions with multiple global minima, addressing a critical need in fields like engineering design, machine learning, and bioinformatics. By incorporating leader-follower dynamics and localized interactions, the algorithm efficiently navigates high-dimensional search spaces to detect multiple optimal solutions. After providing a binary description, a mean-field approximation is derived, and different numerical experiments are conducted to validate the results.
	\end{abstract}


	\section{Introduction}\label{sec:intro}
	Recently, several numerical methods leveraging collective dynamics for minimizing non-convex high-dimensional functions have been developed, \cite{benfenati2022binary,borghi2022consensus,borghi2023constrained,carrillo2018analytical,carrillo2021consensus,Chen2022adamCBO,fornasier2020consensus,fornasier2021consensus,herty2022recent}. These methods, also known as gradient-free methods, offer the advantage of enhancing the efficiency of traditional gradient-based numerical optimization techniques, such as the Newton method or stochastic gradient descent method, \cite{bottou2018optimization}. While these methods can be efficient when gradients are easily computed, they are less practical in scenarios where the function evaluations, which require the gradients computation, are costly or unreliable. In contrast, particle methods explore the search space without needing gradient computations, making them suitable for applications like machine learning and engineering, where simulating physical processes can be expensive. Thus, particle methods provide a more robust approach to find optimal solutions, with an enhanced potential to locate global minima rather than being restricted to local minima.
	Additionally, the introduction of  follower-leader dynamics has been shown to be an effective tool in global optimization problems, allowing a strong speed up of the convergence process, \cite{albi2019leader, albi2023kinetic, motsch2014heterophilious,borghi2023kinetic, albi2021optimized}. Indeed, particle methods are characterized by an exploration phase which slows down the convergence, especially with high diffusion parameter values. By assuming that only a part of the population (followers) is involved in the exploration process, while the remaining one (leaders) quickly reach consensus on the estimated position of the global minimizer, it is possible to achieve faster convergence. 
	With this work, we aim at extending the KBO algorithm enhanced by genetic dynamics (GKBO) introduced in \cite{albi2023kinetic}, to incorporate the capability of exploring objective functions with multiple global minima. Multi-modal optimization is a critical area in numerical optimization where the goal is to identify multiple optimal solutions within an objective function that possesses several global minima. Various approaches, such as the ones described in \cite{preuss2015multimodal}, aim to maintain diversity and prevent premature convergence. Additionally, consensus-based optimization (CBO) methods, like the one introduced in \cite{fornasier2024pde}, offer another effective strategy for handling multiple global minima. This is in contrast to traditional optimization, which typically seeks a single optimal solution. Multi-modal optimization is particularly relevant in fields such as engineering design, machine learning, and bioinformatics, where multiple equally good solutions can exist, each offering unique advantages. A variety of techniques have been proposed for multi-modal optimization, with evolutionary algorithms such as the ones proposed in \cite{goldberg1987genetic, mengshoel2008crowding}. Additionally, polarized CBO methods \cite{bungert2024polarized} have been introduced to further enhance multi-modal optimization. Here, the dynamics is localized using a certain kernel, and particles, instead of being attracted by a common weighted mean, drift towards a weighted mean that emphasizes nearby particles more heavily. In this way, the convergence to a unique global minimizer is prevented, and multiple global minima can be detected. Similarities are shared with models for opinion formation or collective motion, where agents influence strongly their nearest neighbours, leading to the formation of different consensus points or clusters \cite{albi2023topological, during2009boltzmann}.
	The purpose of this paper is to show how to incorporate a similar dynamics in the GKBO algorithm to allow particles to concentrate over the different global minima. The GKBO algorithm belongs to the class of consensus based algorithms, catching ideas from genetic algorithms. The basic assumption is to have a population divided into two groups, namely leaders and followers, similar to the parents and children in the genetic dynamics which mimics biological evolution \cite{zbigniew1993ga,mitchell1995genetic,toledo2014global}. In genetic algorithms, individuals (parents) from the current population are chosen based on their objective values (genetic information) and combined to create the next generation (offspring). This selection process is typically guided by the principle of survival of the fittest, leading the system to progressively approach an optimal solution over successive generations. In this framework, leaders are selected in a manner similar to parents in genetic algorithms, while the generation of children is replaced by followers gradually moving in direction of the leaders' positions. Then, the dynamics of classical KBO is split between the leaders, which move toward the global optimum estimate $\hat x$, and the followers, which explore the minimization landscape \cite{benfenati2022binary,totzeck2021trends}. In particular, each agent is identified by a certain position and label $(x,\lambda)$. If an agent is in the follower status, then $\lambda =0$. On the contrary, if it is in the leaders status, then $\lambda = 1$. In particular two agents with state $(x,\lambda)$ and $(x_*,\lambda_*)$ are supposed to interact as follows
	\begin{equation}\label{eq:GKBO_binary}
		\begin{aligned}
			&x' = x + (\nu_F(x_*-x)+\sigma_FD(x)\xi) \lambda_* (1-\lambda) + \nu_L(\hat{x}(t)-x)\lambda,\\
			&x_*' = x_*,	    
		\end{aligned}
	\end{equation}
	where $x'$ denotes the post interaction position, $\nu_F,\nu_L,\sigma_F$ are positive parameters balancing the exploration and attraction dynamics, $\xi$ is a random perturbation term, $D(x)$ is a diffusion matrix and $\hat{x}(t)$ denotes the estimated position of the global minimizer at time $t$ computed according to the Laplace principle \cite{laplace}. Agents can vary their status in time accordingly to different strategies. In this work we focus on a particular approach, similar to the selection process between parents and children in the genetic algorithm \cite{zbigniew1993ga}. More in details, agents who are in the best position over the cost function are assumed to be leaders, while the others are in the followers status. 
	Consensus-based methods have proven effective for non-convex objectives, but they are limited by their ability to approximate only a single minimizer \cite{fornasier2021global}. To address these limitations, our approach involves establishing distinct clusters, each associated to a certain leader, and to compute distinct estimates of each minimizer across the clusters, based on topological interactions between followers and leaders.
	
	The rest of the paper is organized as follows. In Section \ref{sec:GKBO} we introduce the localised version of the GKBO algorithm, focusing on the description of the rules that govern the dynamics, and of the label switching procedure.  In Section \ref{sec:mean field} we derive the mean field limit, assuming that both the leaders and followers population admit a continuous limit, and in particular the two equations governing the evolution of the followers and leaders densities functions. In Section \ref{sec:num_methods} we introduce the numerical methods showing how to implement a Nanbu type algorithm which mimics the dynamics of the introduced method. In Section \ref{sec:validation}, we perform different simulations, showing the efficiency of the proposed algorithm in terms of success rate and iteration number, and comparing it with cluster-based algorithms.
	
	\section{Localised genetic kinetic based optimization (GKBO) } \label{sec:GKBO} 
	The localized GKBO method catches ideas from the KBO algorithm with genetic dynamics. As in the GKBO algorithm, the population is assumed to be divided into two groups, which are distinguished using labels. However, the dynamics is tailored in such a way that the population divides in different sub-groups clustering around the different global minimizer of the possibly non-convex objective function $\mathcal{E} \colon \RR^d \rightarrow \RR$. Hence, in the long time limit the dynamics solves the global optimization problem given by
	\begin{equation}\label{eq:optProblem}
		\min\limits_{x\in\RR^d} \mathcal{E}(x),
	\end{equation}
	where $\mathcal E$ is assumed to have multiple global minimizer, denoted by $\bar{x}_k$ for any $k=1,\ldots,n_{min}$. In more detail, each agent is described by its position $x \in \mathbb{R}^d$ varying continuously
	and a binary variable for the leadership-level $\lambda \in \{0,1\}$. As in the GKBO, we identify leaders with  $\lambda = 1$ and followers with $\lambda = 0$. We are interested in the
	evolution of the density function \begin{equation}\label{eq:def_f}
		f=f(x,\lambda,t), \qquad f: \RR^d\times\left\{0,1\right\} \times \RR_+\rightarrow \RR_+
	\end{equation}
	where $t\in\RR^+$ denotes the time variable.
	We assume that the total mass of each population is given by
		\begin{equation*}
			\rho_\lambda = \int_{\mathbb{R}^d} f(x,\lambda,t) \, dx,
		\end{equation*}
		for any $\lambda \in \{0, 1\}$, and we suppose $\rho_0 + \rho_1 = 1$. To normalize the population densities, we introduce the scaling factor $1/\rho_\lambda$, and define the normalized density as
		\begin{equation*}
			f_\lambda(x,t) = \frac{f(x,\lambda,t)}{\rho_\lambda}.
		\end{equation*}
		This scaled function $f_\lambda(x,t)$ now satisfies the condition for a probability density function, as it integrates to one over the domain.
		Next, we define the total population density function as
		\begin{equation*}
			g(x,t) = \sum_{\lambda \in \{0,1\}} f(x,\lambda,t),
		\end{equation*}
		which represents the total probability measure for the system.

	\subsection{Binary interaction between agents} 
We introduce $J_c$ clusters of leaders, where each cluster contains a subset of leaders that are grouped based on their proximity. If the number of leaders $N_L$ is finite, then we assume $J_c = N_L$, meaning that each leader forms its own individual cluster. Then, for each cluster $\mathcal{C}_k$, $k = 1,\ldots,J_c$, we define the set of followers associated with the leaders in that cluster. In particular, the set $C_k$ consists of all $x\in \RR^d$ that are closer to the leader in the $k$-th cluster than to any other leader
		\begin{equation}\label{eq:C_k} 
			\mathcal{C}_k = \{x \in \RR^d : \vert x-x_k\vert \leq \vert x-x_j\vert \text{ for any } j = 1,\ldots,N_L, j \neq k\}.
		\end{equation}
		where $x_k$ is the position of the leader in the $k$-th cluster. Thus, each leader is associated with a specific cluster of followers, forming a localized group.
	To avoid ambiguity  and guarantee that each agent belongs to exactly one set $\mathcal{C}_k$, we adopt the convention that if an agent is equidistant from two or more leaders, the agent is assigned to the leader with the lowest index.
	Moreover, we define the estimated position of the global minimizer associated with cluster $k$ at time $t$ as follows
	\begin{equation}\label{eq:x_tot_lead} 
		\hat{c}_k(t) = \frac{\int_{\mathcal{C}_k} x e^{-\alpha \mathcal{E}(x)} g(x,t) \, dx}{\int_{\mathcal{C}_k} e^{-\alpha \mathcal{E}(x)} g(x,t) \, dx},
	\end{equation}
	computed according to Laplace principle \cite{laplace} as a convex combination of particle locations weighted by the cost function. Note that this is a localized average, hence we expect $\hat{c}_k(t) \rightarrow \arg\min\limits_{x\in \mathcal C_k \cap \mathrm{supp }g} \mathcal E$ as $\alpha\rightarrow\infty$. For any $x\in\RR^d$, compute 
	\begin{equation}\label{eq:x_tot} 
		\hat{x}(x;t) = \sum_{k=1}^{J_c} \hat{c}_k(t) \chi(x \in \mathcal{C}_k),
	\end{equation} 
	where $\chi(\cdot)$ is the usual characteristic function. 
	Note that in case we assume $J_c=1$, we have
	\begin{equation}
		\hat{x}(t) = \frac{\int_{\mathbb{R}^d}x e^{-\alpha \mathcal{E}(x)}g(x,t)\,dx}{\int_{\mathbb{R}^d} e^{-\alpha \mathcal{E}(x)}g(x,t)\,dx},
	\end{equation} 
	and Laplace principle yields
	\begin{equation}\label{eq:laplace}
		\lim_{\alpha \to \infty} \Big( -\frac{1}{\alpha} \int_{\mathbb{R}^d} e^{-\alpha \mathcal{E}(x)} g(x,t)\, dx \Big) = \inf_{x\in \text{supp }  g(x,t)} \mathcal{E}(x).
	\end{equation} 
	A binary interaction of agents with state  $(x,\lambda)$ and its nearest neighbor in state $(x_*,\lambda_*)$ is described by their post-interaction positions given by
	\begin{equation}\label{eq:bin_2pop}
		\begin{aligned}
			x' & = x + \Big( \nu_F (x_*-x) + \sigma_F D(x)\xi \Big) \chi_{\mathcal{B}_1(x,x_*;f_1)} \lambda_* (1-\lambda)+ \nu_L(\hat{x}(x;t)-x)\lambda, \\
			x_*' &= x_*,
		\end{aligned}
	\end{equation}
	where $\mathcal{B}_1(x,x_*;f_1)$ is the ball centred in $x$ containing the first nearest neighbor of agent $x$ which is in the leaders status, i.e. agent $x_*$, and $\chi_{\mathcal{B}_1(\cdot)}$ is the usual characteristic function,  $\sigma_F,$ $\nu_F$, $\nu_L$, are positive parameters, $\xi$ is a normally distributed $d$-dimensional random vector, and $D(x)$ is the diffusion matrix, defined to be either 
	\begin{equation}\label{eq:diffusion_iso}
		D(x) = \vert \hat{x}(x;t)-x\vert Id_d,
	\end{equation}
	in the case of isotropic diffusion \cite{pinnau2017consensus}, or 
	\begin{equation}\label{eq:diffusion_ani}
		D(x) = diag\{(\hat{x}(x;t)-x)_1,\ldots,(\hat{x}(x;t)-x)_d\} ,
	\end{equation}
	in the case of anisotropic diffusion \cite{carrillo2021consensus}.
	
	In \eqref{eq:bin_2pop} we use a compact form to describe both the followers and leaders dynamics. In particular, if we assume $\lambda =0$ (and $\lambda_* = 1$) we retrieve the followers dynamics, 
	\begin{align*}
		x' & = x +  \nu_F (x_*-x)\chi_{\mathcal{B}_1(x,x_*;f_1)} + \sigma_F D(x)\xi , \\
		x_*' &= x_*,
	\end{align*}
	where $x,x_*$ denote the pre-interaction positions of a follower and its nearest leader, respectively.
	Followers are attracted toward the leader located in the nearest position and explore the space, searching for a possible position of the global minimizer associated with the selected leader. On the contrary, if $\lambda = 1$, we retrieve the leaders dynamics,
	\begin{align*}
		x' & = x +  \nu_L(\hat{x}(x;t)-x), \nonumber \\
		x_*' & = x_*,
	\end{align*}
	where $x,x_*$ denote the pre-interaction positions of two leaders.
	Hence, leaders do not explore the space but they relax their position toward the estimated position of the global minimizer at time $t$, which is given by $\hat{x}(x;t)$, defined in \eqref{eq:x_tot}. 
	\begin{remark}
			The interaction rule we consider here, is a generalization of binary interactions in the following sense.
			If $\lambda = 0$, the leader $x_*$ influences the post-interaction position of the follower $x$. Specifically, the follower in position $x$ moves towards the leader’s location $x_*$, with the diffusion increasing proportionally to the distance between the two agents. If $\lambda = 1$, the post-interaction position of the leader $x$ is influenced by the pre-interaction positions of the agents in the cluster to which it belongs, which could include the agent $x_*$, contributing to the determination of the consensus point.
	\end{remark}
	
	\section{Derivation of the mean-field equations equations}\label{sec:mean field}
	We assume that each agent modifies its position through a binary interaction occurring with an other agent, who is assumed to be its nearest neighbor in the leaders status. Hence, we introduce the topological ball $\mathit{B}_{r^*}(x,t)$, the ball centred in $x$ whose radius is defined, for a fixed $t \geq 0$, by the following variational problem 
	\begin{equation}\label{eq:radius}
		r^*(x,t) = \arg\min_{\alpha>0} \left\{\int_{\mathit{B}_\alpha (x,t)} f_1(\tilde{x},t) d\tilde{x} \geq \rho^* \ \right\},
	\end{equation}
	where $\rho^*\in(0,1]$ is the target topological mass, namely the ratio $\rho^* = 1/N_L$, where $N_L$ is the total number of leaders, associated to the ball $\mathcal{B}_1(x,x_*;f_1)$ introduced in \eqref{eq:bin_2pop}. 
	For $\lambda \in \{0,1\}$, the evolution of the density function $f_\lambda(x,t)$ is described by the following integro-differential equation
		\begin{equation}\label{eq:boltz_lin}
			\partial_t f_\lambda(x,t)  -\T_\lambda[f](x,t)= \sum_{ \lambda_*\in\{0,1\}} Q(f_\lambda,f_{\lambda_*})(x,t),
		\end{equation}
		where $\T_\lambda[\cdot]$ is the transition operator modeling the change of labels and will be defined in Section \ref{sec:leaders}, and $Q(\cdot,\cdot)$ is the collisional operator. The above equation in weak form, namely for all smooth functions $\phi(x)$ and for any $t\geq 0$ reads 
		\begin{equation}\label{eq:weak_form}
			\begin{aligned}
				\int_{\RR^d} \partial_t f_\lambda(x,t) &\phi(x) \,dx - \int_{\RR^d} \T_\lambda[f](x,t) \phi(x) \, dx = \\ &\eta \sum_{ \lambda_*\in\{0,1\}}   \left \langle \int_{\Omega} \ (\phi(x')-\phi(x)) \, d f_{\lambda_*} \, d f_{\lambda} \right \rangle,
			\end{aligned}
		\end{equation}
		being  $\eta>0$ a constant relaxation rate representing the interaction frequency.
		In order to retrieve the asymptotic behavior of the Boltzmann type equation \eqref{eq:boltz_lin}, we introduce a grazing collision limit argument for the interaction operator, following the approach proposed in \cite{pareschi2013interacting}. The grazing collision limit is connected to the mean-field approximation, especially within the framework of kinetic theory. This limit corresponds to a regime where particle interactions result in very small changes in the agents positions, but at the same time, the number of interactions become very large. As a result, each particle experiences an average effect from the surrounding particles, rather than direct, strong pairwise interactions.  To describe this limit, we introduce a scaling parameter $\varepsilon > 0$, which characterizes the relative strength and frequency of collisions as follows 
		\begin{equation}\label{eq:scaling}
			\begin{split}
				&\nu_F \rightarrow \nu_F \varepsilon,\qquad\nu_L \rightarrow \nu_L \varepsilon, \qquad \sigma_F \rightarrow \sigma_F \sqrt{\varepsilon}, \qquad \eta \rightarrow \frac{1}{\varepsilon}.
			\end{split}
		\end{equation} 
		and we rewrite the scaled binary interactions as 
		\begin{equation}\label{eq:bin_scaled} 
			x'-x = \Big( \nu_F\varepsilon (x_*-x)  + \sigma_F \sqrt{\varepsilon} D(x)\xi\Big) \lambda_* (1-\lambda) + \nu_L \varepsilon (\hat{x}(x;t)-x) \lambda.
		\end{equation}
		As  $\varepsilon\to 0$, the frequency of collisions increases while their individual strength decreases. Then, we introduce a Taylor expansion of the test function $\phi(x)$ in \eqref{eq:weak_form}, as follows
		\begin{equation}
			\phi(x')-\phi(x) =\nabla_x\phi(x)\cdot  (x'-x)  + \frac{1}{2} (x'-x)^T \mathcal{H}_x \phi(x)  (x'-x) + \mathcal{R}(\varepsilon),
		\end{equation} 
		where $\mathcal{H}_x \phi(x)$ denotes the Hessian matrix of the function $\phi(x)$, and $\mathcal{R}(\varepsilon)$ is the remainder term which is given by 
		\begin{equation}\label{eq:remainder}
			\mathcal{R}(\varepsilon) = \frac{1}{2\varepsilon} \sum_{\lambda_*\in \{0,1\}}  \left \langle  \int_{\Omega} (x'-x)^T\mathcal{H}_x  \phi(\bar{x}) (x'-x)  df_{\lambda_*} df_\lambda  \right \rangle , 
		\end{equation}
		with 
		\begin{equation*}
			\bar{x}= \gamma x + (1-\gamma)x',
		\end{equation*}
		for some $\gamma \in [0,1]$. This expansion allows us to capture both the linear and nonlinear effects of the weak, frequent collisions. By doing so, we transition from the microscopic collisional dynamics to an averaged description, which characterizes the mean-field behavior of the system. Then, we rewrite equation \eqref{eq:weak_form} as \begin{equation}\label{eq:weak_formulation_foll}
			\begin{aligned}	
				&\int_{\mathbb{R}^d} \partial_t f_\lambda(x,t)\phi(x) \, dx -\int_{\mathbb{R}^d} \T_\lambda[f](x,t) \phi(x) \, dx =\nonumber \\
				& \sum_{\lambda_*\in\{0,1\}} \left \lbrace  \int_{\Omega} \Big( \nu_F (x_*-x)  \lambda_* (1-\lambda) + \nu_L (\hat{x}(x;t)-x) \lambda \Big)\cdot \nabla_x \phi(x) \, df_{\lambda_*} \, df_{\lambda}\right.\nonumber \\
				&+\left.\frac{\sigma^2_F}{2} \int_{\Omega} [D(x) \mathbf{1}]^T \mathcal{H}_x  \phi(x) [D(x) \mathbf{1}] (1-\lambda)^2\lambda_*^2  \, df_{\lambda_*} \, df_{\lambda} \right \rbrace + \mathcal{R}(\varepsilon),
			\end{aligned}
		\end{equation}
		where $\mathbf{1} = [1,\ldots,1]^T$ is the $d$- dimensional unitary vector, and where we recall that $\Omega =  \RR^d \times \mathcal{B}_{r^*}(x,t)$. By taking  the limit $\varepsilon \to 0$,  we conjecture that one can rigorously demonstrate that $\mathcal{R}(\varepsilon)\to 0$. This is achieved by employing a similar approach to that in \cite{albi2016invisible}, which involves bounding the remainder term with functions of $\varepsilon$. Finally, integration by parts leads to the following equations, which in strong form read
		\begin{equation}
			\begin{aligned}\label{eq:boltz_strong_lead_foll}
				& \partial_t  f_0(x,t) - \T_0[f](x,t)
				= \frac{\sigma_F^2}{2} \Delta_x\Bigl[[D(x) \mathbf{1}]^T [D(x) \mathbf{1}] \rho_1^*(x,t) f_0(x,t) \Bigr] - \\ &\qquad \qquad\qquad \qquad \qquad \qquad \nu_F \nabla_x \cdot \Bigl[\Bigr(m_1^*(x,t) -x \rho_1^*(x,t)\Bigl) f_0(x,t)\Bigr], \nonumber \\&
				\partial_t f_1(x,t) -  \T_1[f](x,t) = -\nu_L   \nabla_x \cdot \Bigl[\Bigr(\hat{x}(x;t) -x \Bigl) \rho_1^*(x,t) f_1(x,t)\Bigr],
			\end{aligned}
		\end{equation}
		where $\Delta_x$ denotes the Laplace operator, and  $D(x)$ is the diffusion matrix defined in \eqref{eq:diffusion_iso}-\eqref{eq:diffusion_ani}, $\mathbf{1} = [1,\ldots,1]^T$ is the $d$-dimensional unitary vector representing the second moment of each component of the vector $\xi$, $\hat{x}(x;t)$ is the global estimate of the global minima at time $t$ defined as in equation \eqref{eq:x_tot}, and
		\begin{equation}\label{eq:mean_leaders}
			m_1^*(x,t) = \int_{\mathcal{B}_{r^*}(x,t)} \tilde{x} \, f_1(\tilde{x},t)\, d\tilde{x}, \qquad
			\rho_1^*(x,t) = \int_{\mathcal{B}_{r^*}(x,t)} f_1(\tilde{x},t) \, d \tilde{x},
		\end{equation}
		denote the centre of mass of the leaders and the leaders density respectively inside the topological ball $\mathcal{B}_{r^*}(x,t)$.
	
	\section{Change of labels}\label{sec:leaders}
	The change of labels is realized with the help of a transition operator which acts as follows
		\begin{equation}\label{eq:master_0}
			\begin{aligned}
				\T_0[f](x,t) &= \pi_{L\to F}(x,\lambda;f)f_1(x,t)-  \pi_{F\to L}(x,\lambda;f) f_0(x,t),  \\
				\T_1[f](x,t) & =  \pi_{F\to L}(x,\lambda;f)f_0(x,t)-  \pi_{L\to F}(x,\lambda;f)f_1(x,t),
			\end{aligned}
	\end{equation}
	where $\pi_{F\to L}(\cdot)$ and $\pi_{L\to F}(\cdot)$ are certain transition rates, possibly depending on the current states. In this context, we consider the weighted strategy proposed in \cite{albi2023kinetic}, similar to the selection process of the genetic algorithm. 
	The idea is to associate a weight $\omega(x,t)$ to each agent in relation to their position $x$ on the function $\mathcal{E}$ that we are going to minimize and to assume that the portion of  agents whose weight is smaller than a certain value $\bar{\omega}$, which depends on the amount of leaders that we would like to generate, is in the leaders status while the remaining ones are in the followers status. 
	Specifically, each agent in position $x$ is associated with a weight  
		\begin{equation*}
			\omega(x,t) =\frac{1}{N} \# \left \lbrace y\in \mathcal{A}(t): \vert  \mathcal{E}(x_{min}) -  \mathcal{E}(y)\vert < \vert  \mathcal{E}(x_{min}) -  \mathcal{E}(x)\vert  \right \rbrace,
		\end{equation*} 
		with
		\begin{equation*}
			x_{min}(t) = \arg \min_{x \in \mathcal{A}(t)} \mathcal{E}(x),
		\end{equation*}
		where $\mathcal{A}(t)$ is the set of agents at time $t$. The value of $\bar{\omega}$ is instead defined to be $N_L/N$, being $N_L$ the total number of leaders that we want to generate.
	The transition rates writes as follows 
	\begin{equation}\label{eq:rates_test_1}
		\begin{aligned}
			\pi_{L\to F}  =\begin{cases}
				1, \quad \text{if } \omega(x,t)>\bar{\omega},\\
				0, \quad \text{if } \omega(x,t)\leq \bar{\omega},\\
			\end{cases}   \qquad
			\pi_{F\to L}  =\begin{cases}
				0, \quad \text{if } \omega(x,t) \geq \bar{\omega},\\
				1, \quad \text{if } \omega(x,t)<\bar{\omega}.\\
			\end{cases}
		\end{aligned}    
	\end{equation}
	The evolution of the emergence and decay of leaders can be described by the master equation
	\begin{equation}\label{eq:Ldef2}
		\dfrac{d}{dt}\rho_\lambda(t) = \int_{\RR^{2d}} \T_\lambda[f](x,t)\, dx,
	\end{equation}
	for $\lambda \in \{0,1\}$, with $\rho_\lambda(t) = \int_{\RR^d} f_\lambda(x,t) dx$, which can be derived starting from \eqref{eq:weak_form}, by taking $\phi(x) = 1$.
	%
	\section{Numerical methods} \label{sec:num_methods}
	In order to solve the mean-field dynamics we consider the Boltzmann-type equation \eqref{eq:boltz_lin} in the scaling limit \eqref{eq:scaling}. This approach allows us to mimic the binary interactions between agents, which is naturally captured by the Boltzmann equation. By adopting this framework, we can describe the evolution of the population densities, incorporating the effects of the interactions. The Boltzmann equation provides a more structured way to model and understand the collective dynamics and emergent behaviors of the population in the mean-field limit. To this aim, we split the evaluation of the dynamics in two different steps, the interaction process and the label evolution, as follows
	\eqref{eq:boltz_lin}
	\begin{align}
		\partial_t f_\lambda(x,t) & = Q^{\varepsilon}_\lambda(f_\lambda,f_\lambda)(x,t),\label{eq:boltz_collision}\\
		\partial_t f_\lambda(x,t) & =\mathcal{T}_\lambda[f](x,t)\label{eq:boltz_label}.
	\end{align}
	In order to approximate the time evolution of the density $f_\lambda(x,t)$ we sample $N_s$ particles $(x_i^0,\lambda_i^0), i=1,\dots,N_s$ from the initial distribution. We consider a time interval $[0, T]$ discretized in $N_t$ intervals of length $h$.
	\paragraph{Interaction step.}
	The interaction step in \eqref{eq:boltz_collision} is solved by means of binary interaction algorithms.
	We denote the approximation of $f_\lambda(x,nh)$ at time $t^n = n h$, by $f_\lambda^{n}(x)$. We decompose the interaction operator introduced in \eqref{eq:boltz_lin} in its gain and loss part,
	\[
	Q^{\varepsilon}_\lambda(f_\lambda,f_\lambda) = \frac{1}{\varepsilon} \Big[Q^{\varepsilon,+}_\lambda(f_\lambda,f_\lambda)- \rho^* f_\lambda \Big],
	\]
	where $\rho^*=1/N_L$ is the topological mass. Considering a forward discretization we obtain
	\begin{equation}\label{eq:collision}
		f_\lambda^{n+1} =\Big( 1-\frac{\rho^*\Delta t}{\varepsilon} \Big)  f_\lambda^n + \frac{\rho^*\Delta t}{\varepsilon}  \frac{Q^{\varepsilon,+}(f_\lambda^n,f_\lambda^n )}{\rho^*}.
	\end{equation} 
	Equation \eqref{eq:collision} can be interpreted as follows. With probability $1-\rho^*\Delta t/\varepsilon$ an individual in position $x$, with label $\lambda$ does not interact with other individuals and, with probability $\rho^*\Delta t /\varepsilon$, it interacts with another individual according to 
	\begin{equation}\label{eq:bin_2pop_disc}
		\begin{aligned}
			x^{n+1} = x^n + &\Big(\varepsilon \nu_F (x_*^n-x^n) + \sqrt{\varepsilon} \sigma_F D(x^n)\xi \Big)\lambda_*^n (1-\lambda^n)+ \\ & +
			\varepsilon \nu_L(\hat{x}(x^n,t^n)-x^n)\lambda^n,
		\end{aligned}
	\end{equation}
	where $x_*^n$ is selected to be the nearest neighbour to agent $x$ in the leader status. We assume $\rho^*\Delta t = \varepsilon$ to maximize the total number of interactions and ensure that at each time step all agents interact with another individual with probability one.
	\paragraph{Labels evolution.}
	In order to simulate changes of the label $\lambda$, we discretize equation \eqref{eq:boltz_label}. For any fixed $x\in \RR^{d}$, we obtain
	\begin{align}\label{eq:lambda_evolution}
		f_0^{n+1}(x) = (1-\varepsilon ~\pi_{F\to L})~ f_0^n(x) + \varepsilon ~\pi_{L\to F} ~f_1^n(x),\nonumber\\
		f_1^{n+1}(x) = (1-\varepsilon ~\pi_{L\to F})~ f_1^n(x) + \varepsilon ~\pi_{F\to L}~ f_0^n(x),
	\end{align}
	where $\pi_{F\to L}(\cdot)$ and $\pi_{L\to F}(\cdot)$ are the transition rates as defined in \eqref{eq:rates_test_1}. 
	
	The details of the numerical scheme are summarized in Algorithm \ref{alg_binary}. Here, the parameters $\delta_{stall}$ and $j_{stall}$ are used to check if consensus has been reached in the last $j_{stall}$ iterations within a tolerance $\delta_{stall}$. In more detail, we stop the iteration if the distance of the current and previous mean $\hat{x}(x;t)$, for any $x$ is smaller than the tolerance $\delta_{stall}$ for at least $j_{stall}$ iterations. 
	\begin{alg}~[Localized GKBO]\label{alg_binary}
		\begin{enumerate}
			\item[\texttt 1.] Draw $(x_i^0,\lambda_i^0)_{i=1,\dots,N_s}$ from the initial distribution $f^0_\lambda(x)$ and set $n=0$, $j=0$, $j_i=0$, for any $i=1,\ldots,N_s$.
			\item[\texttt 2.] \texttt{while} $n<N_t$ \texttt{and} $j<j_{stall}$ 
			\begin{enumerate}
				\item   Associate to each leader its group of followers by computing the set $\mathcal{C}_k$ as in \eqref{eq:C_k}.
				\item \texttt{for} $i=1$ \texttt{to}  $N_s$
				\begin{itemize}
					\item For each follower, select the nearest leader with position $y^{n}_k$, $k\neq i$.
					\item Compute the estimated position of the global minima $\hat{c}_k^n$ as in \eqref{eq:x_tot_lead}, for any $k=1,\ldots,N_L$, being $N_L$ the total number of leaders to be generated.  
					\item Associate to each agent $i$ a certain $\hat{x}(x_i^n;t^n)$ as in \eqref{eq:x_tot}. 
					\item Compute the new positions as in \eqref{eq:bin_2pop_disc}.
					\item Compute the following probabilities rates 
					\[
					p_{L} =\varepsilon ~\pi_{F\to L}(x_i^{n+1},\lambda_i^n), \qquad 	p_{F}=\varepsilon ~\pi_{L\to F}(x_i^{n+1},\lambda_i^n),
					\]  
					with $\pi_{F\to L}, \pi_{L\to F}$ as in \eqref{eq:rates_test_1}.
					\begin{itemize}
						\item \texttt{if} $\lambda_i^n = 0$,\\  with probability $p_{L}$ agents $i$ becomes a leader: $\lambda_i^{n+1} = 1$.
						\item \texttt{if} $\lambda_i^n = 1$,\\  with probability $p_{F}$ agents $i$ becomes a follower: $\lambda_i^{n+1} = 0$.
					\end{itemize}
				\end{itemize}
				\item Compute $\hat{x}(x_i^{n+1};t^{n+1})$ as in equation \eqref{eq:x_tot}.
				\item  \texttt{if} $\lVert \hat{x}(x_i^{n+1};t^{n+1})-\hat{x}(x_i^{n};t^{n})\rVert_\infty\leq \delta_{stall}$
				\begin{enumerate}
					\item [] $j_i\leftarrow j_i+1$
					\item [] \texttt{end if}
				\end{enumerate}
				
				\texttt{end for}
				\item [] $j = \min_i\{j_i\}$,
				\item [] $n\leftarrow n+1$.
			\end{enumerate}
			\texttt{end while} 
		\end{enumerate}
	\end{alg}
	The above algorithm is inspired from Nanbu's method \cite{nanbu1980direct}, for larger class of direct simulation Monte-Carlo algorithm for interacting particle dynamics we refer to \cite{albi2013binary,pareschi2013interacting}. 
	
	\section{Validation tests}\label{sec:validation} 
	In this section we test the performance of the localized GKBO algorithm in terms of number of detected minima, success rate and number of needed iterations. 
	We run $M=20$ simulations and, according to \cite{benfenati2022binary,carrillo2018analytical}, we consider a simulation successful if for any minimizer $\bar{x}$ there exists a weighted mean $\hat{x}(x;t)$ s.t.
	\begin{equation}\label{eq:success}
		\lVert \hat{x}(x;t) - \bar{x}\rVert_{\infty} \leq 0.25.
	\end{equation}
	We set $\alpha = 5\cdot 10^6$ and we adopt the numerical log-sum-exp trick described in \cite{fornasier2021consensus} to allow for arbitrary large values of $\alpha$. We assume $N_s=600$ and that agents are initially uniformly distributed in the hypercube $[-10,10]^d$.  At time $t=0$ we suppose all agents are in the followers status and that labels change according to the weighted criterion defined in Section \ref{sec:leaders}, with $\bar{\omega} = N_L/N$, being $N_L$ the total number of leaders that we would like to generate, which varies in the different tests. We let the dynamics in \eqref{eq:bin_2pop} to evolve for $N_t=10000$ iterations with $\varepsilon=0.1$. We set $j_{stall} = 1000$, $\delta_{stall}= 10^{-4}$.  We assume $\nu_F=1$, $\nu_L= 2$ while the diffusion parameter and the dimension change in the different tests and will be specified later.
	In the following tests we construct multi-modal versions of Rastrigin and Ackley function as benchmark functions, \cite{jamil2013literature}, and we consider anisotropic diffusion.
	
\begin{figure}[H]
		\centering
		\includegraphics[width=5.88cm]{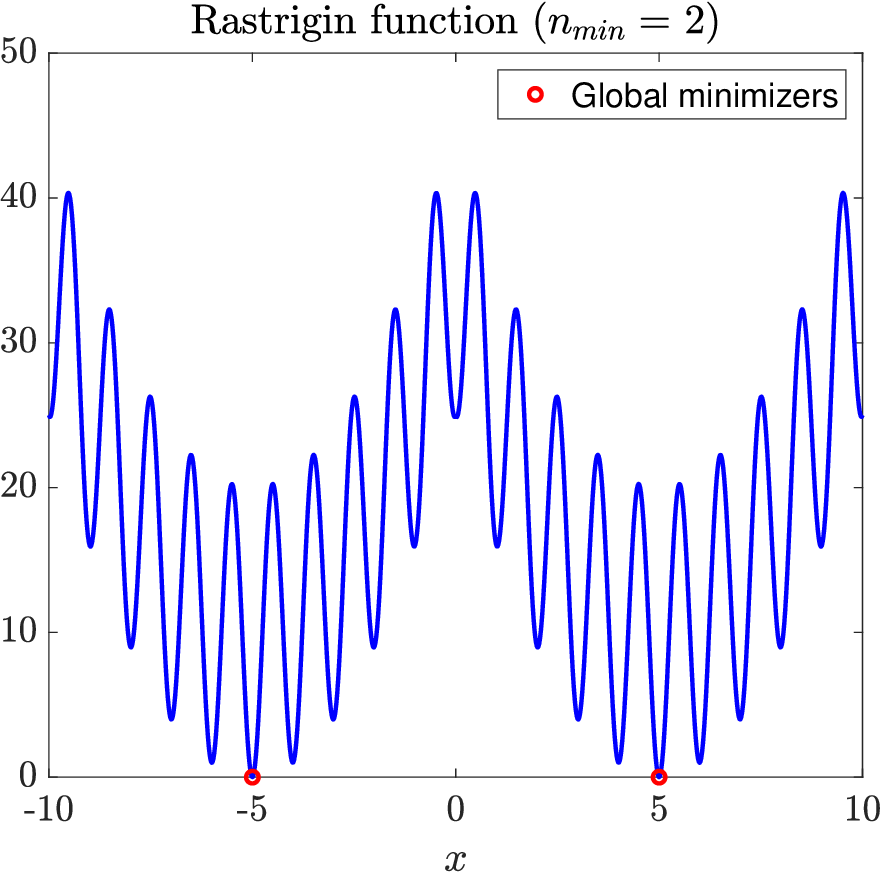}
		\qquad
		\includegraphics[width=5.88cm]{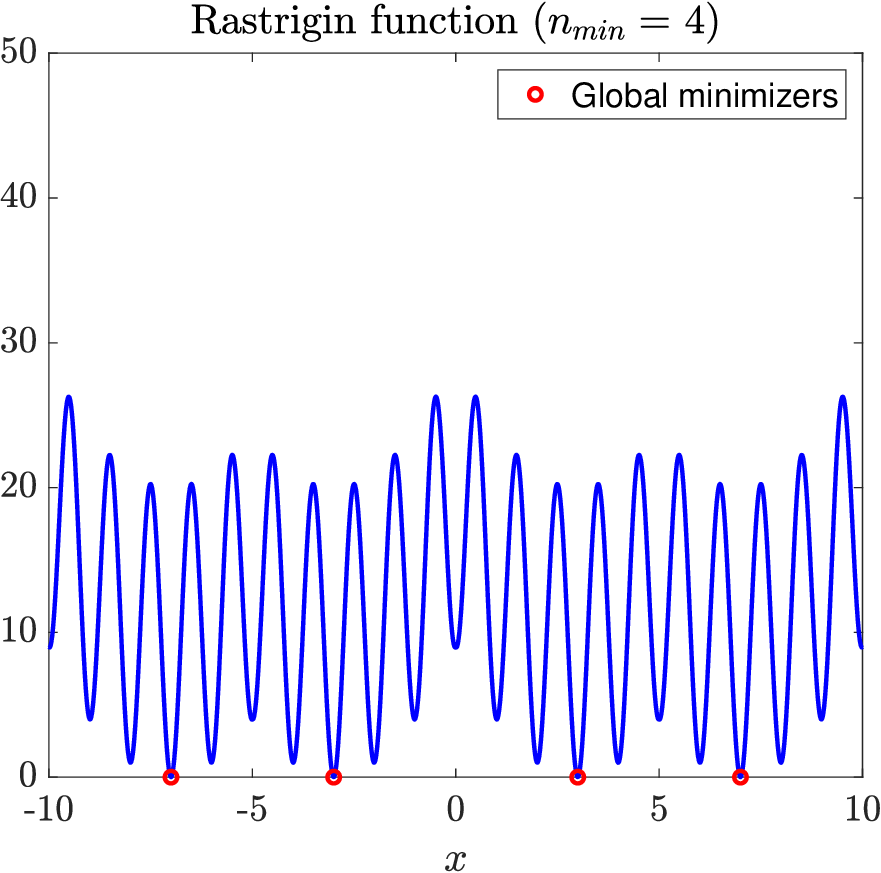}
		\caption{Multi-modal Rastrigin function \eqref{eq:multi_rastrigin}. On the left, with two global minima and on the right with four global minima. We highlight with markers their position. }\label{fig:rastrigin}
	\end{figure}

	\paragraph{Multi-modal Rastrigin function.}
	To construct a multi-modal Rastrigin function, we start by considering the classical Rastrigin function, which has a single global minimizer at $\bar{x} = 0$, expressed as follows
	\begin{equation}\label{eq:rastrigin}
		R(\xx) = \frac{1}{d}\sum_{i=1}^d (x_i^2 - 10 \cos(2\pi x_i)),
	\end{equation}
	where $d$ is the dimension, and $\xx = (x_1,\ldots,x_d)$. Then we define its multi-modal version as follows
	\begin{equation}\label{eq:multi_rastrigin}
		R_M(\xx) = \min_{\xx \in \RR^d} \{R(\xx-\bar{x}_1),\ldots,R(\xx-\bar{x}_{n_{min}})\}, 
	\end{equation}
	being $n_{min}$ the number of global minima $\bar{x}_k$, $k=1,\ldots,n_{min}$ that we would like to introduce. In Figure \ref{fig:rastrigin} on the left the multi-modal Rastrigin function with $n_{min} = 2$ global minima in $\bar{x}_1 = -5$, $\bar{x}_2 = 5$, and on the right the same function but with $n_{min} = 4$ global minima in $\bar{x}_1 = -7$, $\bar{x}_2 = -3$, $\bar{x}_3=3$, $\bar{x}_4=7$. We highlight with markers the position of each $\bar{x}_k$, $k=1,\ldots,n_{min}$. 
	\paragraph{Multi-modal Ackley function.}
	Similar to the Rastrigin function, we consider the uni-modal Ackley function
	\begin{equation}
		A(\xx) = -20 \exp\Bigg( -0.2\sqrt{\frac{\sum_{i=1}^d x_i^2}{d}} \Bigg) -\exp\Bigg( \frac{\sum_{i=1}^d \cos{(2\pi x_i)}}{d}\Bigg)  + 20+e,
	\end{equation} 
	which has a unique global minimizer in $\bar{x} = 0$,
	where $d$ is the dimension of the problem and $\xx = (x_1,\ldots,x_d)$. Its multi-modal version reads as follows
	\begin{equation}\label{eq:multi_ackley}
		A_M(\xx) = \min_{\xx \in \RR^d} \{A(\xx-\bar{x}_1),\ldots,A(\xx-\bar{x}_{n_{min}})\}, 
	\end{equation}
	being $n_{min}$ the number of global minima $\bar{x}_k$, $k=1,\ldots,n_{min}$ that we would like to introduce. In Figure \ref{fig:rastrigin} on the left the multi-modal Ackley function with $n_{min} = 2$ global minima in $\bar{x}_1 = -3$, $\bar{x}_2 = 3$, and on the right the same function but with $n_{min} = 4$ global minima in $\bar{x}_1 = -7$, $\bar{x}_2 = -3$, $\bar{x}_3=3$, $\bar{x}_4=7$. We highlight with markers the position of each $\bar{x}_k$, $k=1,\ldots,n_{min}$. 
	
	\begin{figure}[h!]
		\centering 
		\includegraphics[width=5.88cm]{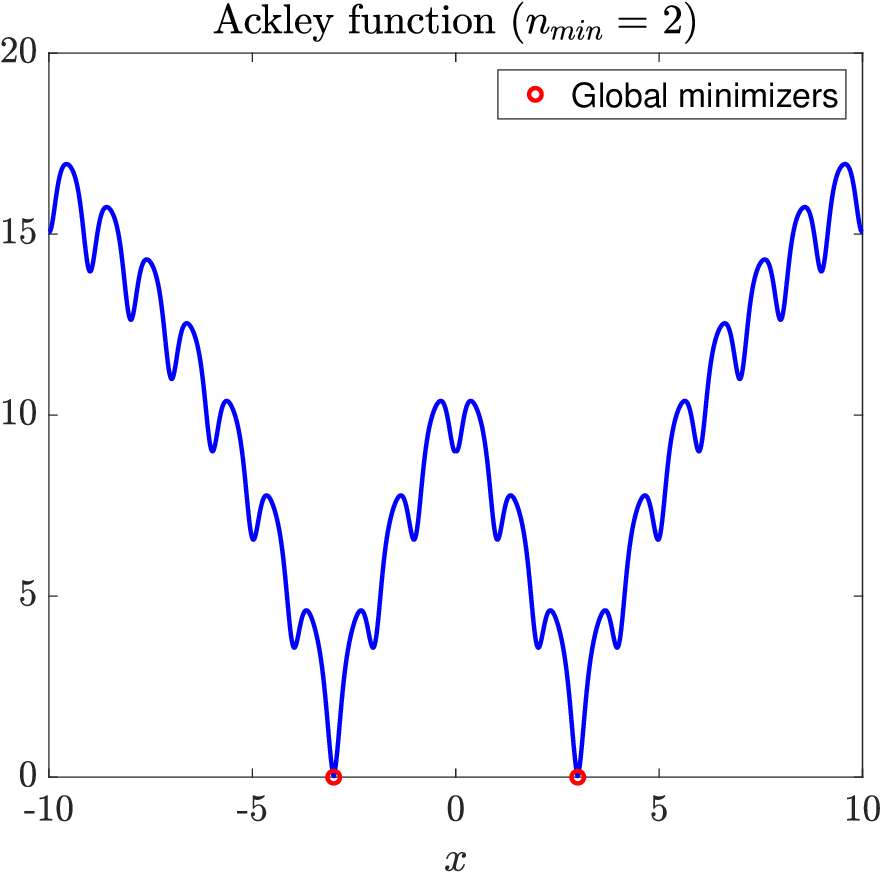}
		\qquad
		\includegraphics[width=5.88cm]{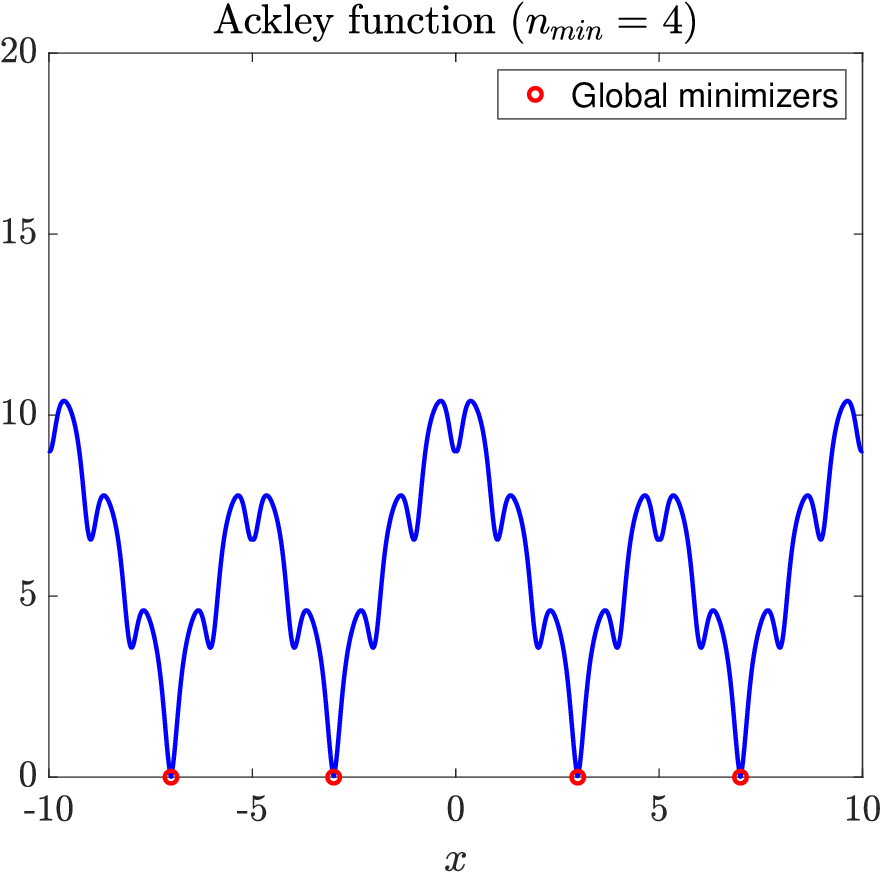}
		\caption{Multi-modal Ackley function \eqref{eq:multi_ackley}. On the left, with two global minima and on the right with four global minima. We highlight with markers their position. }\label{fig:ackley}
	\end{figure}
	
	\subsection{Test 1: Number of detected minima and leaders dependence}\label{sec:test1} 
\begin{figure}[h!]
		\centering 
		\includegraphics[width=.459\textwidth]{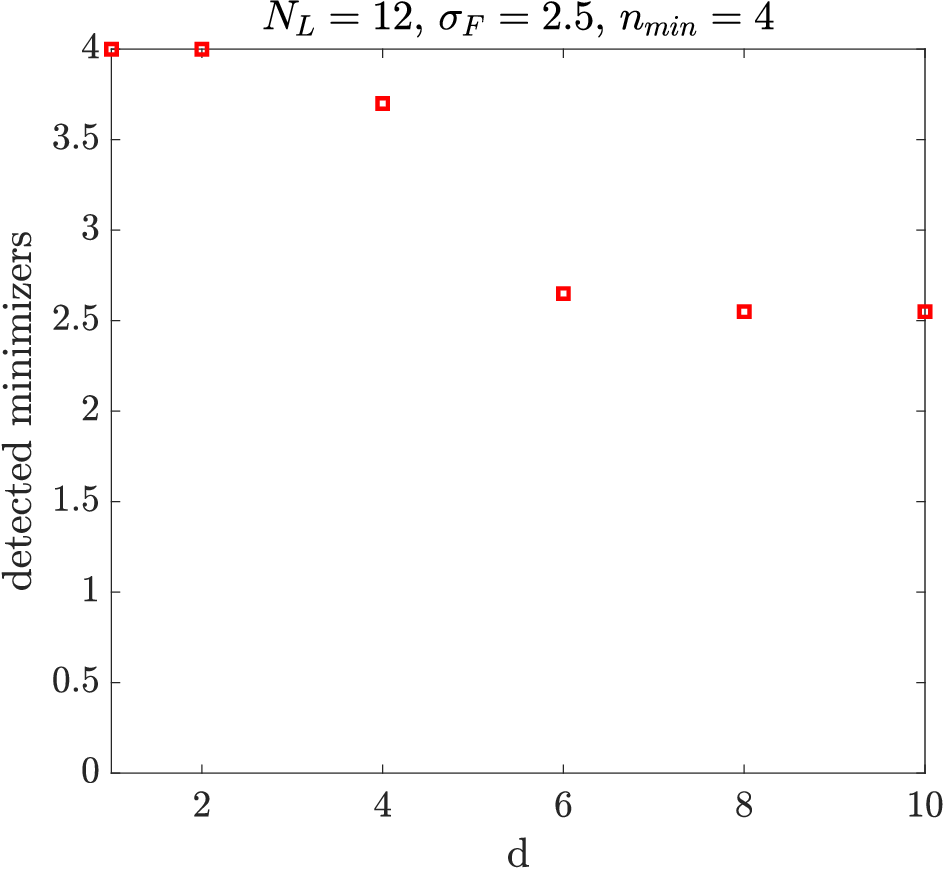}
		\caption{Average number of detected minima as the dimension $d$ varies for $\sigma_F= 2.5$, for the Rastrigin function \eqref{eq:multi_rastrigin} with $n_{min} = 4$ global minima. The markers denote the average number of detected minima for different values of $d$. }
		\label{fig:num_det_lead}
	\end{figure}

\begin{figure}[h!]
		\centering 
		\includegraphics[width=5.88cm]{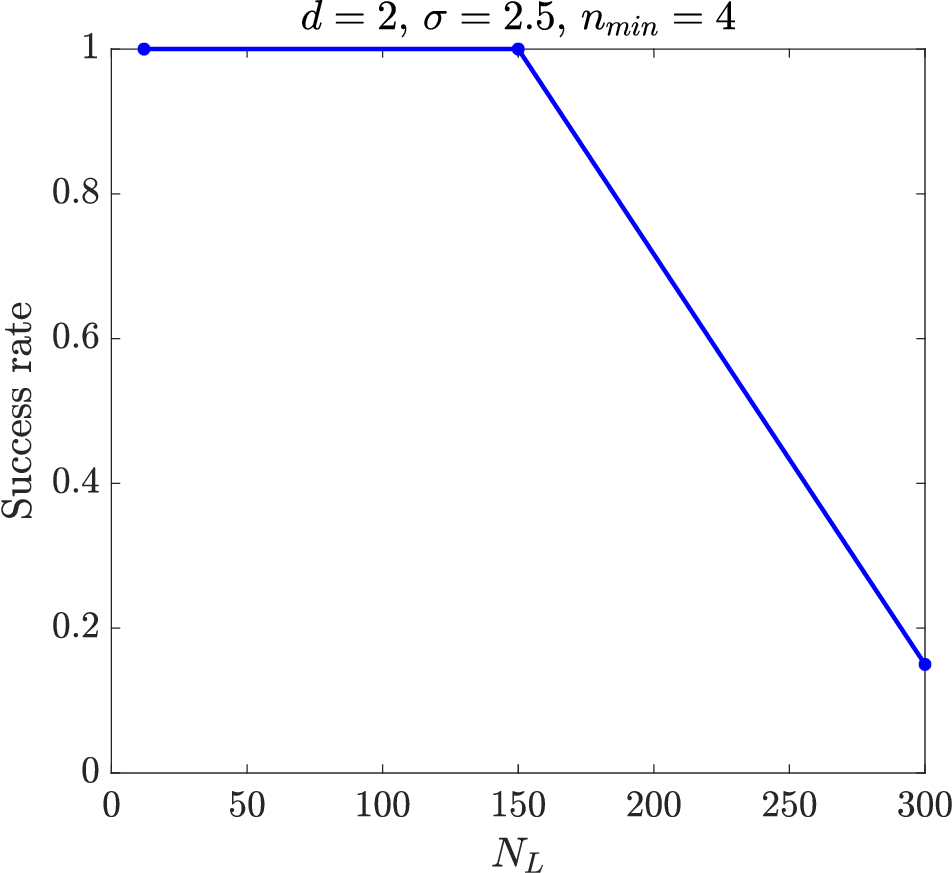}
		\qquad
		\includegraphics[width=5.88cm]{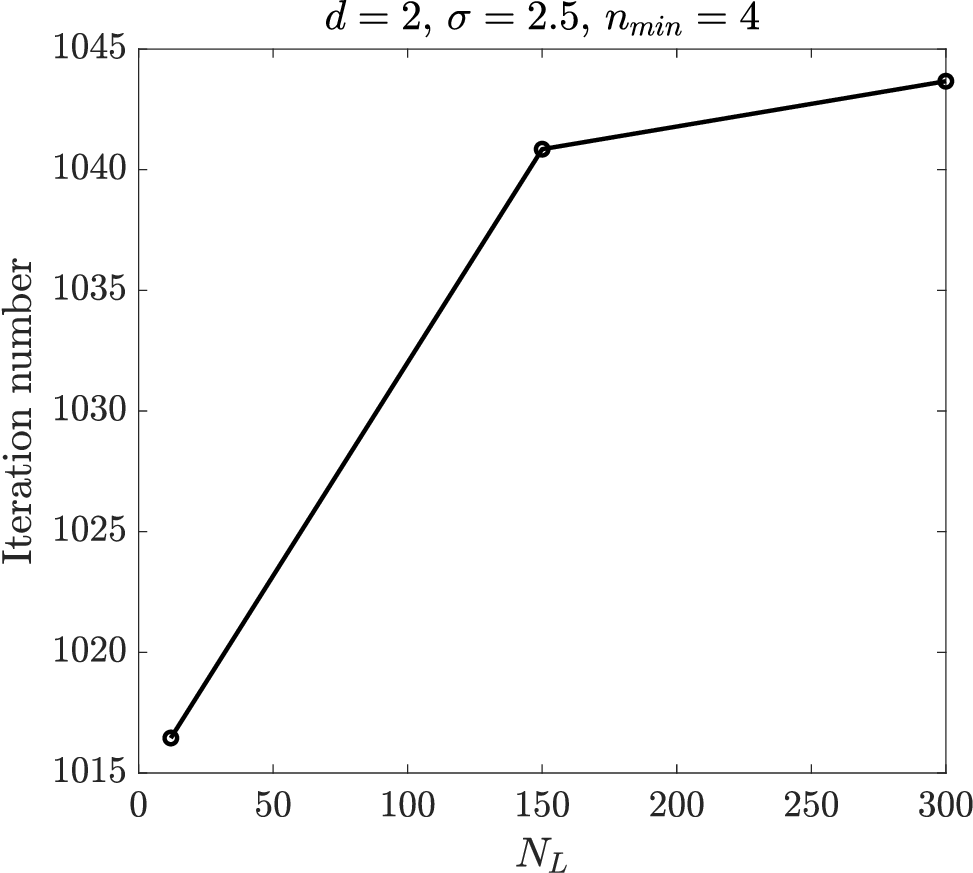}
		\caption{Success rate and number of iterations as $N_L$ varies for the Rastrigin function \eqref{eq:multi_rastrigin} with $n_{min} =4$ minima for $d=2$ and $\sigma_F=2.5$. The markers denote the value of the success rates and mean number of iterations. }
		\label{fig:NL_range}
	\end{figure}

	We first consider the Rastrigin function with $n_{min} =4$ global minima. We test the effectiveness of the localized GKBO algorithm in detecting the whole set of global minima as the dimension $d$ varies. In Figure \ref{fig:num_det_lead} we report the average number of detected minima as a function of $d$. By increasing the dimension, the probability of detecting the total amount of minima decreases since the complexity of the problem increases.
	Secondly, we fix the dimension $d=2$, and the diffusion parameter $\sigma_F= 2.5$ and we test the effectiveness of the GKBO algorithm in detecting $n_{min} = 4$ minima as the total number of leaders varies. We choose $N_L = 12$, $N_L = 120$, $N_L = 300$.
	In Figure \ref{fig:NL_range} on the left the success rate and on the right the iteration number as $N_L$ varies. As the number of leaders increases, the success rate decreases while the iteration number remains almost the same. A high number of leaders does not imply better and/or faster convergence. 

	
	\subsection{Test 2: Comparison for different dimension and diffusion parameters } \label{sec:test2}
\begin{figure}[h!]
		\centering 
		\includegraphics[width=5.88cm]{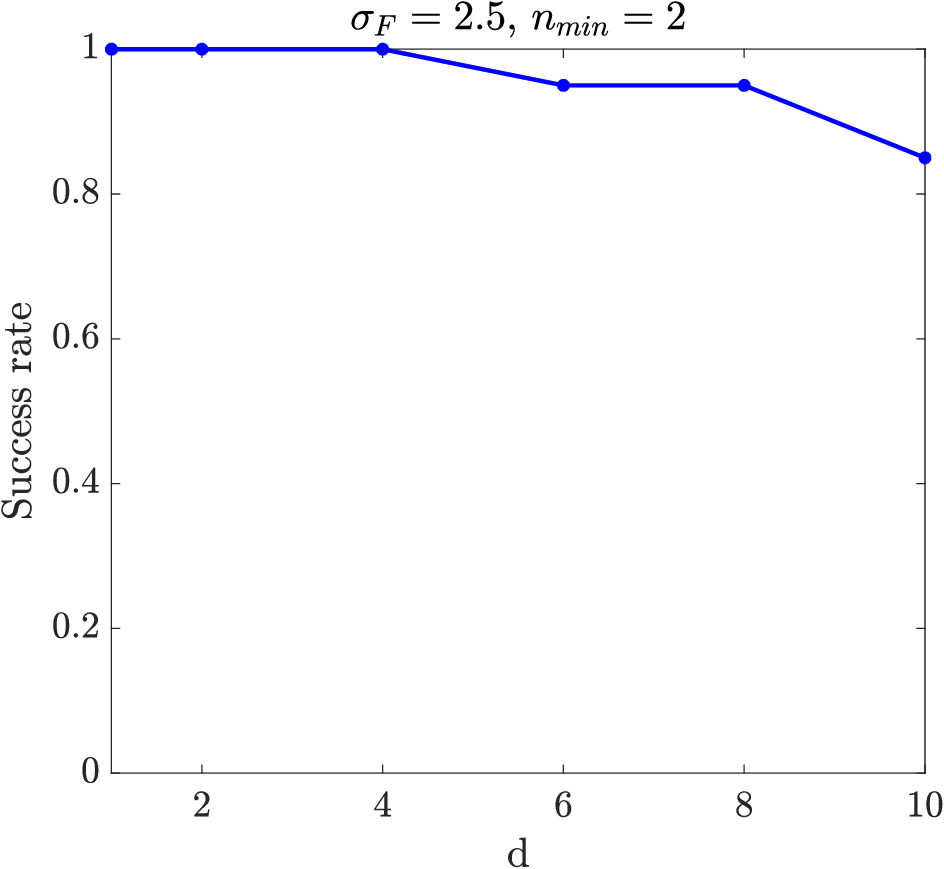}
		\qquad
		\includegraphics[width=5.88cm]{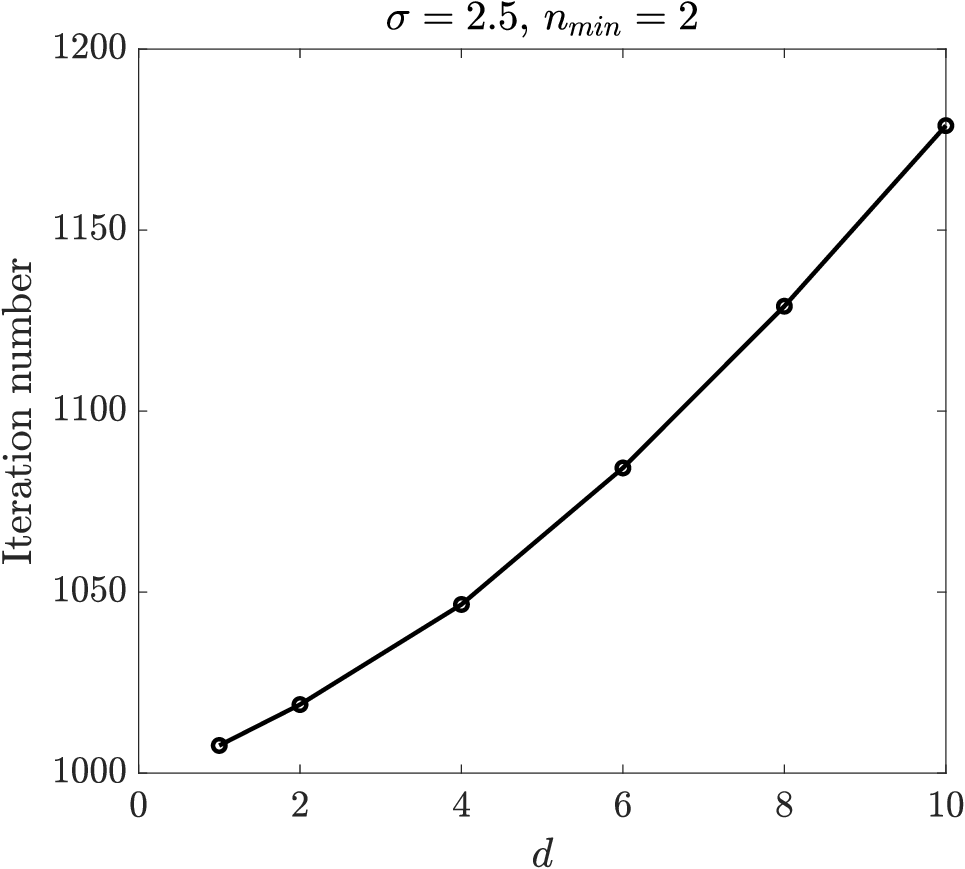}
		\caption{Success rate and number of iterations as $d$ varies for the Rastrigin function \eqref{eq:multi_rastrigin} with $n_{min} =2$ minima and $\sigma_F=2.5$. The markers denote the value of the success rates and mean number of iterations.}
		\label{fig:d_range}
	\end{figure}

	We first fix $\sigma_F=2.5$ and let the dimension vary as $d = 1,\ldots,10$. 
	In Figure \ref{fig:d_range} the success rates and average number of iterations obtained with the localized GKBO in the case of the Rastrigin function with $n_{min} = 2$ global minima is shown. The success rate and iteration number increase proportionally to the dimension. 
	
\begin{figure}[h!]
		\centering 
		\includegraphics[width=5.88cm]{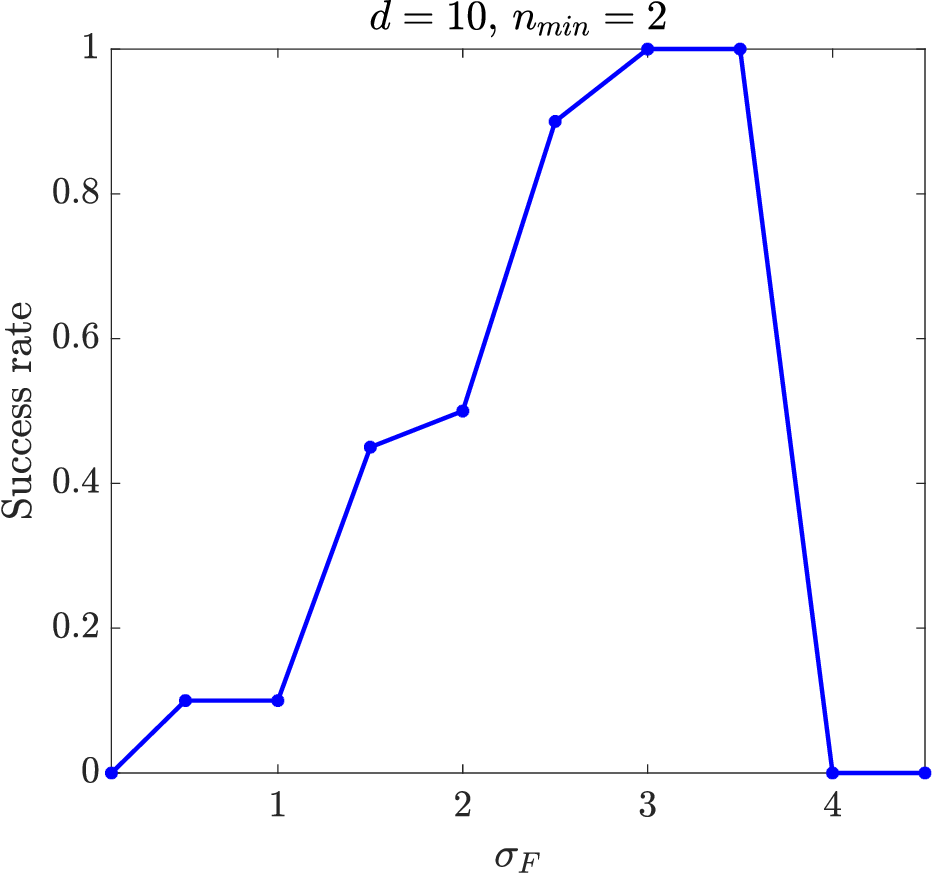}
		\qquad
		\includegraphics[width=5.88cm]{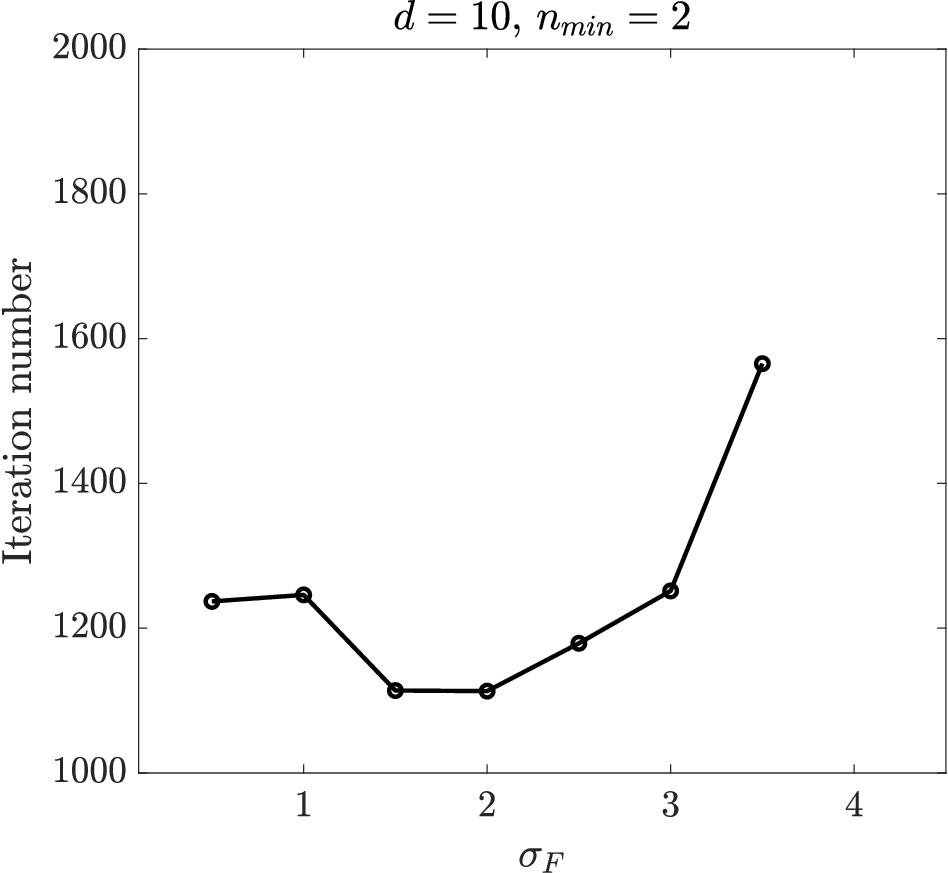}
		\caption{Success rate and number of iterations as $\sigma_F$ varies for the Rastrigin function \eqref{eq:multi_rastrigin} with $n_{min} =2$ minima and $d = 10$. The markers denote the value of the success rates and mean number of iterations.  
		}
		\label{fig:sigma_range}
	\end{figure}

	Secondly, we fix $d=10$ and let the diffusion parameter to vary as $\sigma_F = 0.1,\ldots,4.5$. 
	In Figure \ref{fig:sigma_range} the success rates and average number of iterations obtained with the localized GKBO in the case of the Rastrigin function with $n_{min} = 2$ global minima is shown. The success rate and iteration number is strictly dependent of the choice of the diffusion parameter. A small diffusion parameter causes particles to become trapped in local minima, while a large diffusion parameter prevents the formation of consensus, increasing the iteration number.
	
	\subsection{Test 3: Comparison between the localized GKBO and the polarized CBO algorithms} \label{sec:test4}
\begin{figure}[h!]
		\centering 
		\includegraphics[width=5.88cm]{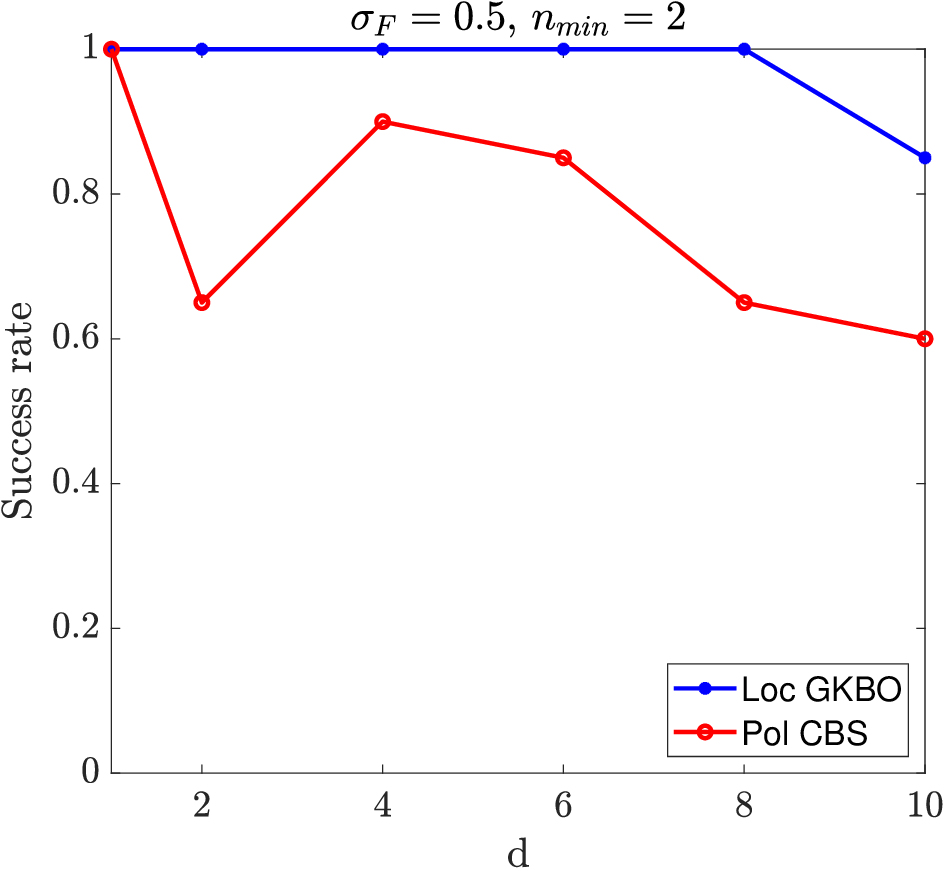}
		\qquad
		\includegraphics[width=5.88cm]{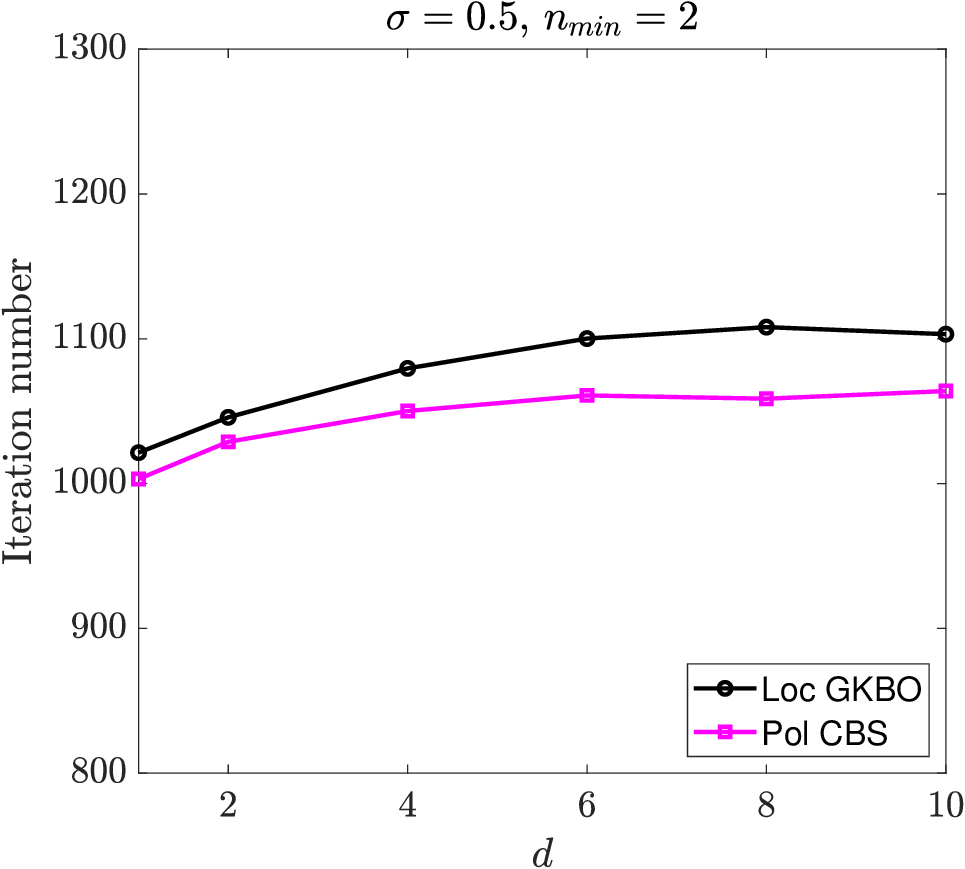}
		\caption{Comparison between the localized GKBO \eqref{eq:bin_2pop} and the polarized CBO algorithms \eqref{eq:CBS} in terms of success rate and number of iterations as $d$ varies for the Ackley function \eqref{eq:multi_ackley} with $n_{min} =2$ minima and $\sigma_F=0.5$. The markers denote the value of the success rates and mean number of iterations.}
		\label{fig:comparison}
	\end{figure}

	In this Section we compare the polarized CBO introduced in \cite{bungert2024polarized}, and the localized GKBO. The idea of the polarized CBO is to localize the dynamics as follows. Assume to have one population, and introduce $J_c$ clusters whose position at time $t$ is given by 
	\begin{equation}
		c_j(t) = \frac{\sum_{i=1}^N p_{ij} x_i e^{-\alpha \mathcal{E}(x_i)}}{\sum_{i=1}^N p_{ij} e^{-\alpha \mathcal{E}(x_i)}},
	\end{equation}
	where $p_{ij}$ denotes the probability of particle $i$ to belong to cluster $j$ for any $i=1,\ldots,N$, $j=1,\ldots,J_c$. Then the estimated position of the global minima at time $t$ is given by 
	\begin{equation}
		\hat{x}_i(t) = \sum_{j=1}^{J_c} p_{ij} c_j(t), 
	\end{equation}  
	and each particle $x_i$ for any $i=1,\ldots,N$ modifies its position as follows 
	\begin{equation}\label{eq:CBS}
		x_i' = x_i + \nu (\hat{x}_i(t) - x_i) + \sigma D(x_i)\xi, 
	\end{equation}
	where $\xi$ is a normaly distributed random number, $D(x)$ is the usual diffusion matrix, and $\nu,\sigma$ are positive parameters. At time $t=0$, the probability $p_{ij}$ is drawn from a uniform distribution, that is $p_{ij}\sim \mathcal{U}([0,1])$. Then at each time, we assume 
	\begin{equation}
		p_{ij} = \begin{cases}
			1, \qquad \text{if } |x_i-x_j| < |x_i-x_k|, \quad \text{for any } k \neq j, k = 1,\ldots, J_c,\\
			0, \qquad \text{else.} 
		\end{cases}
	\end{equation}
	Note that the localized GKBO resembles the polarized CBO approach. However, the inclusion of a follower-leader dynamics leads to an enhancement in terms of efficiency, as shown in Figure \ref{fig:comparison}, where the two approaches are compared. In particular, we consider the multi-modal Ackley function with $n_{min} = 2$ global minima and we set $\sigma_F= \sigma = 0.5$, $\nu_F=\nu=1$, $N_L=J_c=4$ and test the results as the dimension varies between $d=1,\ldots,10$. The same test has been conducted on the multi-modal Rastrigin function, showing similar results. 
	
	\section{Conclusion}\label{sec:conclusion}
	In this study, we have extended the KBO algorithm enhanced with genetic dynamics (GKBO) to its localized version, useful to detect multiple global minima.  Our approach, which divides the population into leaders and followers, shows improvements in solving complex global optimization problems, especially when compared to the polarized CBO method, as demonstrated in the numerical experiments.
	The validation tests indicate that the localized GKBO algorithm effectively detects global minima across varying dimensions and diffusion parameters. Specifically, the success rates and iteration numbers scale proportionaly with the problem's dimensionality, indicating the algorithm's robustness and adaptability in the tested dimensions. Moreover, the variation in diffusion parameters further illustrates the algorithm's flexibility and efficiency in diverse optimization landscapes.
	Overall, the localized GKBO algorithm represents an advancement in global optimization techniques, combining the strengths of kinetic and genetic algorithms to achieve a good performance in detecting multiple global minima.
	
	
	\section*{Acknowledges}
		This work has been written within the activities of GNCS group of INdAM (Italian National Institute of High Mathematics).
		FF was partially supported by  the Italian Ministry of University and Research (MUR) through the PRIN 2020 project (No. 2020JLWP23) “Integrated Mathematical Approaches to Socio–Epidemiological Dynamics".

	
	\bibliographystyle{emss}

\begin{thebibliography}{99}
		
		\bibitem{albi2016invisible}
		G. Albi, M. Bongini, E. Cristiani, D. Kalise.
		\newblock Invisible control of self-organizing agents leaving unknown environments.
		\newblock In {\em SIAM Journal on Applied Mathematics}, 76, No. 4, 1683-1710 ( 2016).
		
		\bibitem{albi2019leader}
		G. Albi, M. Bongini, F. Rossi, and F. Solombrino.
		\newblock Leader formation with mean-field birth and death models.
		\newblock In {\em Math. Models Methods Appl. Sci.},29, No. 4, 633-679 (2019).
		
		
		\bibitem{albi2023topological}
		G. Albi, and F. Ferrarese.
		\newblock Kinetic description of swarming dynamics with topological interaction and transient leaders.
		\newblock {arXiv preprint \em arXiv:2307.12044} (2024).
		
		\bibitem{albi2023kinetic}
		G. Albi, F. Ferrarese, and C. Totzeck.
		\newblock Kinetic-based optimization enhanced by genetic dynamics.
		\newblock In {\em Math. Models Methods Appl. Sci. },  33, No. 14, 2905-2933 (2023). 
		
		\bibitem{albi2021optimized}
		G. Albi, F. Ferrarese and C. Segala.
		\newblock Optimized leaders strategies for crowd evacuation in unknown environments with multiple exits.
		\newblock In {\em Crowd Dynamics, Volume 3: Modeling and Social Applications in the Time of COVID-19}, pp. 97-131 (2021).
		
		\bibitem{albi2013binary}
		G. Albi and L. Pareschi.
		\newblock Binary interaction algorithms for the simulation of flocking and
		swarming dynamics.
		\newblock In {\em Multiscale Model. Simul.}, 11, No. 1, 1-29 (2013).
		
		
		\bibitem{benfenati2022binary}
		A. Benfenati, G. Borghi, and L. Pareschi.
		\newblock Binary interaction methods for high dimensional global optimization
		and machine learning.
		\newblock In {\em Appl. Math. Optim.}, 86, No. 1, 41p. (2022).
		
		
		\bibitem{borghi2022consensus}
		G. Borghi, M. Herty, and L. Pareschi.
		\newblock A consensus-based algorithm for multi-objective optimization and its
		mean-field description.
		\newblock  arXiv preprint { \em  arXiv:2203.16384}(2022).
		
		\bibitem{borghi2023constrained}
		G. Borghi, M. Herty, and L. Pareschi.
		\newblock Constrained consensus-based optimization.
		\newblock In {\em SIAM J. Optim.}, 33, No. 1, 211-236 (2023).
		
		\bibitem{borghi2023kinetic}
		G. Borghi, and L. Pareschi.
		\newblock Kinetic description and convergence analysis of genetic algorithms for global optimization.
		\newblock arXiv preprint { \em arXiv:2310.08562}, (2023).
		
		\bibitem{bottou2018optimization}
		L. Bottou, F.~E. Curtis, and J. Nocedal.
		\newblock Optimization methods for large-scale machine learning.
		\newblock In {\em SIAM Rev. 60}, No. 2, 223-311 (2018).
		
		\bibitem{bungert2024polarized}
		L. Bungert, T. Roith and W. Philipp.
		\newblock Polarized consensus-based dynamics for optimization and sampling.
		\newblock In {\em Math. Program.}, pp 1-31, (2024).
		
		\bibitem{carrillo2018analytical}
		J.~A. Carrillo, Y.~P. Choi, C. Totzeck, and O. Tse.
		\newblock An analytical framework for consensus-based global optimization
		method.
		\newblock In {\em Math. Models Methods Appl. Sci.},
		28, No. 6, 1037-1066 (2018).
		
		
		\bibitem{carrillo2021consensus}
		J.~A. Carrillo, S. Jin, L. Li, and Y. Zhu.
		\newblock A consensus-based global optimization method for high dimensional
		machine learning problems.
		\newblock In {\em ESAIM, Control Optim. Calc. Var. }, 27, Suppl., Paper No. S5, 22 p. (2021).
		
		
		\bibitem{Chen2022adamCBO}
		S. Chen, J. Jin and L. Lyu.
		\newblock A consensus-based global optimization method with adaptive momentum
		estimation.
		\newblock In {\em Commun. Comput. Phys. }, 31, No. 4, 1296-1316 (2022)
		
		\bibitem{laplace}
		A. Dembo. 
		\newblock Large deviations techniques and applications.
		\newblock In {\em Springer} (2009).
		
		
		\bibitem{during2009boltzmann}
		B. D{\"u}ring, P. Markowich, and M.~T. Wolfram.
		\newblock Boltzmann and Fokker--Planck equations modelling opinion formation in the presence of strong leaders.
		\newblock In {\em Proc. R. Soc. Lond., Ser. A, Math. Phys. Eng. Sci.} No. 2112, 3687-3708 (2009).
		
		\bibitem{fornasier2020consensus}
		M. Fornasier, H. Huang, L. Pareschi, and P. S{\"u}nnen.
		\newblock Consensus-based optimization on hypersurfaces: Well-posedness and
		mean-field limit.
		\newblock In {\em Math. Models Methods Appl. Sci.},
		30, No. 14, 2725-2751 (2020).
		
		\bibitem{fornasier2021consensus}
		M. Fornasier, H. Huang, L. Pareschi, and P. S{\"u}nnen.
		\newblock Consensus-based optimization on the sphere: Convergence to global
		minima and machine learning.
		\newblock In {\em J. Mach. Learn. Res.}, 22, No. 237, 1--55 (2021).
		
		\bibitem{fornasier2021global}
		M. Fornasier, T. Klock, and K. Rield.
		\newblock Consensus-based optimization methods converge globally.
		\newblock {arXiv preprint \em arXiv:2103.15130}, (2021).
		
		\bibitem{fornasier2024pde}
		M. Fornasier, and L. Sun.
		\newblock A PDE Framework of Consensus-Based Optimization for Objectives with Multiple Global minima.
		\newblock{arXiv preprint \em arXiv:2403.06662}, (2024).
		
		\bibitem{goldberg1987genetic}
		D.E. Goldberg, and J. Richardson.
		\newblock Genetic algorithms with sharing for multimodal function optimization.
		\newblock In { \em Genetic algorithms and their applications: Proceedings of the Second International Conference on Genetic Algorithms}, Vol. 4149 (1987).
		
		\bibitem{herty2022recent}
		M. Herty, E. Iacomini, and G. Visconti.
		\newblock Recent trends on nonlinear filtering for inverse problems.
		\newblock In {\em  Commun. Appl. Ind. Math. }, 13, No. 1, 10-20 (2022).
		
		\bibitem{jamil2013literature}
		M. Jamil and X.~-S. Yang.
		\newblock A literature survey of benchmark functions for global optimisation
		problems.
		\newblock In {\em Int. J. Math. Model. Numer. Optim.}, 4, No. 2, 150-194 (2013).
		
		
		\bibitem{mengshoel2008crowding}
		O.J. Mengshoel, and D.E. Goldberg.
		\newblock The crowding approach to niching in genetic algorithms.
		\newblock In {\em  Evolutionary computation}, 16, No. 3, 315-354 (2008).
		
		\bibitem{zbigniew1993ga}
		Z. Michalewicz.
		\newblock  Genetic Algorithms + Data Structures = Evolution Programs.
		\newblock In {\em Springer}, 1996.
		
		\bibitem{mitchell1995genetic}
		M. Mitchell.
		\newblock Genetic algorithms: An overview.
		\newblock In {\em Complexity}, 1, No. 1, 31-39 (1995).
		
		\bibitem{motsch2014heterophilious}
		S. Motsch and E. Tadmor.
		\newblock Heterophilious dynamics enhances consensus.
		\newblock In {\em SIAM Rev.},  56, No. 4, 577-621 (2014).
		
		\bibitem{nanbu1980direct}
		K. Nanbu.
		\newblock Direct simulation scheme derived from the {B}oltzmann equation. {I}. {M}onocomponent gases
		\newblock In { \em J. Phys. Soc. Jpn.}, 49 No. 5, 2042-2049, (1980).
		
		\bibitem{pareschi2013interacting}
		L. Pareschi and G. Toscani.
		\newblock  Interacting multiagent systems: kinetic equations and {M}onte {C}arlo methods.
		\newblock {\em OUP Oxford}, (2013).
		
		\bibitem{pinnau2017consensus}
		R. Pinnau, C. Totzeck, O. Tse, and S. Martin.
		\newblock  A consensus-based model for global optimization and its mean-field
		limit.
		\newblock In {\em Math. Models Methods Appl. Sci.},
		27 No. 1, 183-204, (2017).
		
		\bibitem{preuss2015multimodal}
		M. Preuss.
		\newblock Multimodal optimization by means of evolutionary algorithms.
		\newblock{\em Springer}, (2015)
		
		\bibitem{toledo2014global}
		C.~F.~M. Toledo, L. Oliveira, and P.~M. Fran{\c{c}}a.
		\newblock Global optimization using a genetic algorithm with hierarchically
		structured population.
		\newblock In {\em J. Comput. Appl. Math.}, 261, 341-351 (2014).
		
		\bibitem{totzeck2021trends}
		C. Totzeck.
		\newblock Trends in consensus-based optimization.
		\newblock In {\em  Active particles. Volume 3. Advances in theory, models, and applications.}, pages 201--226. Springer, (2021).
		
	\end{thebibliography}

\end{document}